\documentclass [12 pt, letterpaper] {article}

        \usepackage {times,amssymb,amsmath,array,calc,enumerate}
   {\begin {list}{} %{\textbf{Example}\arabic{example}.}
   {\usecounter{example}%
   \small%
   \setlength{\labelsep}{0pt}\setlength{\leftmargin}{40pt}%
   \setlength{\labelwidth}{0pt}%
   \setlength{\listparindent}{0pt}}}%
   {\end{list}}

\renewcommand {\theequation} {\@arabic\c@equation}

\makeatother
   {\begin {list}{} %{\textbf{Example}\arabic{example}.}
   {%\usecounter{example}
   \small%
    \setlength{\itemindent}{10pt}
    \setlength{\labelsep}{5pt}\setlength{\leftmargin}{40pt}%
    \setlength{\labelwidth}{5pt}%
   \setlength{\itemsep}{-1.8pt}
    \setlength{\listparindent}{0pt}
   }
   }%
   {\end{list}}

\newcommand{\B}{\mathcal{B}}

\newcommand{\A}{\mathcal{A}}
\newcommand{\M}{\mathcal{M}}

\newcommand{\pb}{\vspace{0.3cm}\noindent }
\newcommand{\N}{\mathcal{N}}

\newcommand{\R}{\mathcal {R}}

\renewcommand {\theequation} {\arabic{equation}}
\newtheorem{Def}{Definition}[section]
\newtheorem{Th}{Theorem}[section]
\newtheorem{Col}{Corollary}[section]
\newtheorem{Prop}{Proposition}[section]
\newtheorem{lemma}{Lemma}[section]
\begin {document}

\vspace{2cm}

\centerline{\Large Maximal Injective Subalgebras of Tensor Products
}

 \vspace{0.3cm}

 \centerline{\Large of Free Group Factors }

\vspace{2cm}

\centerline{\large Junhao Shen}

\vspace{1cm}

%\centerline{\large Dedicated to Richard. V. Kadison on the occasion
%of his eightiesth birthday}

%\vspace{1cm}

\centerline{ Department of Mathematics, University of New Hampshire,
Durham NH, 03824.}

\vspace{0.3cm}

\centerline{email: jog2@cisunix.unh.edu}

\vspace{2cm}

 \noindent{\bf Abstract:} In this article, we proved the following results.
Let $L(F(n_i))$ ($n_i\ge 2$) be the free group factor on $n_i$
generators and $\lambda (g_{i})$ be one of standard generators of
$L(F(n_i))$ for $1\le i\le N$. Let $\A_i$ be the abelian von Neumann
subalgebra of $L(F(n_i))$ generated by $\lambda(g_{i})$. Then the
abelian von Neumann subalgebra $\otimes_{i=1}^N\A_i$ is a maximal
injective von Neumann subalgebra of $\otimes_{i=1}^N L(F(n_i))$.
When $N$ is equal to infinity, we obtained McDuff factors that
contain maximal injective abelian von Neumann subalgebras.

\thispagestyle{empty} %\vfill
%\newpage

\vspace{2cm}

\section{Introduction}

Let $\mathcal H$ be a separable complex Hilbert space, $\mathcal
B(\mathcal H)$ be the algebra of all bounded linear operators from
$\mathcal H$ to $\mathcal H$. A von Neumann algebra $\R$ is called
``injective" if it is the range of a norm one projection from
$\mathcal B(\mathcal H)$ onto $\R$. Since injective von Neumann
algebras form a monotone class, it follows that any injective von
Neumann algebra is contained in some maximal injective von Neumann
algebra $\mathcal R_1$.

In his influential list of problems presented at the conference in
Baton Rouge in 1967, R. Kadison asked ([3]; Problem 7) whether each
self-adjoint operator in a II$_1$ factor lies in some hyperfinite
subfactor. This problem was answered in the negative in a remarkable
paper [9] by S. Popa. More specifically, if $L(F(n))$ is the free
group factor on   $n$ generators, and $\lambda(g)$ is one of the
standard generators of $L(F(n))$, then S. Popa showed that the
abelian von Neumann subalgebra generated by $\lambda(g)$ is a
maximal injective von Neumann subalgebra of $L(F(n))$. It follows
that $\lambda(g)$ is not contained in any hyperfinite subfactor of
$L(F(n))$, which solves Kadison's problem as mentioned  above.
%It was expected by S. Popa that every non-atomic finite injective
%von Neumann algebra is $*-$isomorphic to a maximal injective
%subalgebra of each nonhyperfinite type II$_1$ factor.
 In [2], L. Ge
provided more examples of maximal injective von Neumann subalgebras
in type II$_1$ factors. Actually he showed that each non-atomic
injective finite von Neumann algebra with a separable predual is
maximal injective in its free product with any von Neumann algebra
associated with a countable discrete group. %Combining with K.
%Dykema's result in [2], it follows that every non-atomic finite
%injective von Neumann algebra is $*-$isomorphic to a maximal
%injective subalgebra of each free group factor.
Note that the type II$_1$ factors listed in [9], [2] that contain
maximal injective abelian subalgebras are all non-$\Gamma$ factors.

 In this paper, we provide   examples of maximal injective von
Neumann subalgebras in   McDuff factors of type II$_1$. In
particular, we consider maximal injective von Neumann subalgebras in
tensor products of free group factors. By developing the techniques
from [9], [2], we are able to prove the following result.

\pb {\bf Theorem } {\em  Suppose  $\{n_i\}_{i=1}^N$ is a sequence of
integers where $n_i\ge 2$ for all $1\le i\le N$ and $N$ is finite or
infinite. Let $F(n_i)$ be the free group with the standard
generators $\{g_{i,j}\}_{j=1}^{n_i}$ for   $1\le i\le N$. Let the
group $G $ be $\times_{i=1}^N F(n_i)$, the direct product of all
$F(n_i)$'s. And $F(n_i)$ is identified with its canonical image in
$G$. Let $\lambda$ be the left regular representation of $G$ and $\M
(=L(G)\cong \otimes_{i=1}^NL(F(n_i)))$ be the group von Neumann
algebra associated with $G$. Let $\A$ be the abelian von Neumann
subalgebra of $\M$ generated by the unitary elements
$\{\lambda(g_{i,1})| 1\le i\le N\}$.  Then $\A$ is a maximal
injective subalgebra of $\M$, thus  not contained in any hyperfinite
subfactor of $\M$.}

\vspace{0.3cm}

When $N$ is equal to infinity, we obtain  examples of McDuff factors
of type II$_1$  (for example, $\otimes_{i=1}^\infty L(F(2))$) that
contain maximal injective abelian von Neumann subalgebras.
% Combining
%with L. Ge's result in [3], we shows that that every non-atomic
%finite injective von Neumann algebra is $*-$isomorphic to a maximal
%injective subalgebra of some McDuff factor of type II$_1$.

The organization of the paper is as follows. We introduce some basic
knowledge in section 2. One useful lemma by R. Kadison is quoted. In
 section 3, some technical lemmas   needed in later section are proved.
  In section 4, we prove our main theorem, Theorem 4.1, of the
paper.

\vspace{0.3cm} It was expected by S. Popa that every non-atomic
finite injective von Neumann algebra is $*-$isomorphic to a maximal
injective subalgebra of each nonhyperfinite type II$_1$ factor. We
hope that our work will provide some new insights into S. Popa's
question.

\vspace{0.3cm} The author wishes to express his deep gratitude to
Prof. L. Ge for many stimulating and fruitful conversations.

\vspace{2cm}

\section{Preliminaries}

Let $\mathcal H$ be a separable complex Hilbert space, $\mathcal
B(\mathcal H)$ be the algebra of all bounded liner operators from
$\mathcal H$ to $\mathcal H$. (For the general theory of operator
algebras, we refer to [5] and [8].) A von Neumann algebra $\R$ is
called ``injective" if it is the range of a norm one projection from
$\mathcal B(\mathcal H)$ onto $\R$. Since injective von Neumann
subalgebras of a von Neumann algebra $\M$ form a monotone class, it
follows that any injective von Neumann subalgebra of $\M$ is
contained in some maximal injective von Neumann subalgebra $\mathcal
R_1$ of $\M$.

\vspace{0.3cm}

Let $\M$ be a finite von Neumann algebra with the tracial state
$\tau$. If $\omega$ is a free filter on $\Bbb N$ then denote by
$\mathcal M^{\omega}$ the quotient of the von Neumann algebra
$l^{\infty}(\Bbb N, \mathcal M)$ by the $0$-ideal of the trace
$\tau_{\omega}$, $ \tau_{\omega} ((x_n)_n)=\lim_{n\rightarrow
\omega}\tau(x_n)$. Then $\mathcal M^{\omega}$ is a finite von
Neumann algebra, $\tau_{\omega}$ is a trace on $\mathcal
M^{\omega}$. $\M$ is naturally embedded in $\M^\omega$ as the
algebra of constant sequence. (see [10])

\vspace{0.3cm} A separable type II$_1$ factor $\M$ has the property
$\Gamma$ of Murray and von Neumann (see [7]) if for any $x_1,
\ldots, x_n\in \M$, $\epsilon>0$ there exists a unitary element
$u\in \M$ such that $\tau(u)=0$ and $||ux_i-x_iu||_2\le \epsilon, \
1\le i \le n.$

\vspace{0.3cm} It is well-known that a separable type II$_1$ factor
$\M$ has the property $\Gamma$ of Murray and von Neumann if and only
if $\M'\cap \M^\omega\ne \Bbb CI$. If $\M\cong \M\otimes \R_0$, then
$\M$ is called a McDuff factor, where $\R_0$ is the unique
hyperfinite type II$_1$ factor. It is known in [6] that $\M$ is a
McDuff factor if and only if $\M'\cap \M^\omega$ is non-commutative.
Since $\R_0\otimes \R_0\cong \R_0$, $\R_0'\cap \R_0^\omega$ is
noncommutative.

\vspace{0.3cm}

Let $\R$ be an injective von Neumann subalgebra of   $\M$. Then
$\mathcal R$ can be decomposed as $ \mathcal R_1 \oplus \mathcal
R_2$, where $R_1$ is a type I von Neumann subalgebra of $\mathcal M$
and $\mathcal R_2$ is a     type II$_1$ von Neumann subalgebra of
$\mathcal M$. From Connes's celebrated result (see [1]), both of
$\mathcal R_1$ and $\mathcal R_2$ are injective. Then we have
\begin{lemma}
  If $\mathcal R_2 \ne 0$, then $\mathcal R'\cap \mathcal
  R_2^{\omega}$ $(\subset \mathcal R'\cap \mathcal M^{\omega})$ is
  noncommutative.
\end{lemma}
\noindent {\bf Proof: } From Lemma XVI.1.5 of [11], it follows that
 $\mathcal R_2 \cong \mathcal Z \otimes \mathcal R_0$,
where $\mathcal Z$ is the center of $\mathcal R_2$ and $\mathcal
R_0$ is the        hyperfinite factor of type II$_1$.  % It is easy to see
%that $I_{\Z} \otimes (\R_0'\cap \R_0^\omega) \ \subset \ \mathcal
%R_2'\cap \mathcal
%  R_2^{\omega}\  \subset \
%\mathcal R'\cap \mathcal
 % R_2^{\omega}$.
    From the fact that  $\R_0'\cap \R_0^\omega$ is a  non-commutative von Neumann algebra, we obtain that $\mathcal R'\cap \mathcal
  R_2^{\omega}$ ($\subset \mathcal R'\cap \mathcal M^{\omega}$) is
  also
  noncommutative.

% (see
%Takesaki, ``Theorey of Operator Algebras III Lemma XVI.1.5) \qquad Q.E.D.

\vspace{0.3cm}

As a corollary, we have
\begin{Col} Let $\M$ be a finite von Neumann algebra with the tracial state
$\tau$. Let $\R$ is an injective von Neumann subalgebra of $\M$. Let
$\A$ be an abelian self-adjoint von Neumann subalgebra of $\R$.
Suppose $\mathcal R= \mathcal R_1 \oplus \mathcal R_2$, where $\R_1$
is a type I von Neumann subalgebra of $\mathcal M$ and $\mathcal
R_2$ is a type II$_1$ von Neumann subalgebra of $\mathcal M$. If
$\R_2 \ne 0$, then there exists an element $x\in \mathcal R'\cap
\mathcal R_2^{\omega}$ not contained in $\mathcal A^{\omega}$.
\end{Col}

\vspace{0.3cm}

Here, we quote a useful lemma from [4].

\begin{lemma} {\em (From [4])}
If $\mathcal N_0$ is a countably decomposable von Neumann algebra,
$\mathcal N$ is the von Neumann algebra of $n\times n$ matrices with
entries in $\mathcal N_0$, and $\mathcal S$ is an abelian
self-adjoint subset of $\mathcal N$, then there is a unitary element
(matrix) $u$ in $\mathcal N$ such that $uau^{-1}$ has all its
non-zero entries on the diagonal for each $a$ in $\mathcal S$.
\end{lemma}

\begin{Def} Let $\M$ be a finite von Neumann algebra with the tracial
state $\tau$. Let $\A$ be a von Neumann subalgebra of $\M$. Then the
normalizer, $N(\A)$, is defined as the set consisting of all unitary
elements $u$ in $\M$ such that $u\A u^*=\A$. If $\M =\{N(\A)\}''$,
then $\A$ is called to be regular in $\M$. If $\A=\{N(\A)\}''$, then
$\A$ is called to be singular in $\M$.
\end{Def}

\vspace{0.3cm} The following lemma tells us that a maximal abelian
von Neumann subalgebra $\A$ in a finite type I von Neumann algebra
$\M$ has to be  regular in $\M$.

\begin{lemma}
Suppose $\mathcal M$ is a  finite   type I von Neumann algebra and
$\mathcal A$ is a maximal abelian  von Neumann subalgebra of
$\mathcal M$. Then
 $\mathcal A$ is regular in $\mathcal M$, i.e.
$\mathcal M$ is generated by the normalizer  of $\mathcal A$.
\end{lemma}
\noindent {\bf Proof: }Decompose $\mathcal M$ as $\oplus_i \mathcal
Z_i\otimes M_{n_i}(\Bbb C)$, where $\mathcal Z_i$ is   abelian
  von Neumann subalgebra of $\M$. It is sufficient to
show the following statement: suppose $\mathcal A_i$ is a maximal
abelian von Neumann subalgebra in $\mathcal Z_i\otimes M_{n_i}(\Bbb
C)$, then $\mathcal Z_i\otimes M_{n_i}(\Bbb C)$ is generated by the
normalizer of $\mathcal A_i$. Since $\mathcal A_i$ is an abelian
  von Neumann subalgebra, $\mathcal A_i$ is generated by
a self-adjoint element $x$ in $\mathcal A_i$. By Lemma 2.2 we know
that there exists a unitary element $u$ in $\mathcal Z_i\otimes
M_{n_i}(\Bbb C)$ such that $uxu^*$ is a diagonal matrix. i.e. $uxu^*
= diag (x_1, x_2, \ldots, x_n)$ where $x_j$ is in $\mathcal Z_i$.
Since $u\mathcal A_iu^*$ is also a maximal abelian   von Neumann
subalgebra in $\mathcal Z_i\otimes M_{n_i}(\Bbb C)$ and generated by
$uxu^* = diag (x_1, x_2, \ldots, x_n)$, we easily have that
$u\mathcal A_iu^*= \mathcal Z_i \otimes \mathcal D_i$ where
$\mathcal D_i$ is the von Neumann subalgebra generated by
$\{e_{ss}\}_{1\le s\le n_i}$ in $M_{n_i}(\Bbb C)$ and
$\{e_{st}\}_{1\le s,t\le n_i}$ is the canonical system of matrix
units in $M_{n_i}(\Bbb C)$. It follows that $\mathcal Z_i\otimes
M_{n_i}(\Bbb C)$ is generated by the normalizer of $u\mathcal
A_iu^*$; consequently by the normalizer of $ \mathcal A_i $. \qquad
Q.E.D.

\vspace{0.3cm}

\begin{lemma} We have the following statements.
\begin{enumerate}[(1)]
\item Suppose $\mathcal R_2$ is a       type II$_1$  injective von Neumann
algebra  with the traical state $\tau$.
 Suppose $\mathcal A$ is a  maximal abelian   von Neumann subalgebra
 of $\mathcal R_2$. Then there exists a unitary element $w$ in
 $\mathcal R_2$ such that $w$ is orthogonal to $\mathcal A$ in $\mathcal R_2$ with respect to
 $\tau$, i.e. $\tau(wx)=0$ for all $x\in \A$.
\item Let $\mathcal R$ be $ \mathcal R_1\oplus \mathcal
R_2$ with the tracial state $\tau$, where $ \mathcal R_1$ is a type
I von Neumann algebra and $ \mathcal R_2$ is a nonzero
 type II$_1$  injective von Neumann algebra. Let $\mathcal A$ is a maximal
abelian   von Neumann subalgebra of $\mathcal R$.  Let $P_1, P_2$ be
the central supports of $\mathcal R_1, \mathcal R_2$, respectively.
Then there exists a unitary $w$ in $\mathcal R_2$,  so a partial
isometry in $\mathcal R$ with $ww^*=w^*w=P_2$, such that $w$ is
orthogonal to $\mathcal A$ in $\mathcal R$ with respect to $\tau$.
\end{enumerate}
\end{lemma}

\noindent {\bf Proof: } We need only to show (1), since (2) follows
from (1) directly. Note that $\mathcal R_2$ can be decomposed as
$\mathcal Z \otimes \mathcal R_0$ where $\mathcal Z$ is the center
of $\mathcal R_2$ and $\mathcal R_0$ is the unique hyperfinite
II$_1$ factor. It is easy to see that $\mathcal R_0 \cong \mathcal
R_0\otimes M_2(\Bbb C)$. Thus $\mathcal R_2$ can also be viewed as
$(\mathcal Z \otimes \mathcal R_0) \otimes M_2(\Bbb C)$ and
$\mathcal A$ is a maximal abelian   von Neumann subalgebra in
$(\mathcal Z \otimes \mathcal R_0) \otimes M_2(\Bbb C)$. By Lemma
2.3, there exists a unitary element $u$ in $(\mathcal Z \otimes
\mathcal R_0) \otimes M_2(\Bbb C)$ such that $uau^*$ has all its
non-zero entries on the diagonal for each $a$ in $\mathcal A$.
Suppose $\{ e_{ij}\}_{1\le i,j\le 2}$ is the canonical system of
matrix units of $M_2(\Bbb C)$ in $(\mathcal Z \otimes \mathcal R_0)
\otimes M_2(\Bbb C)$. Let $v =I_{\mathcal Z \otimes \mathcal
R_0}\otimes e_{12}+ I_{\mathcal Z \otimes \mathcal R_0}\otimes
e_{21}$ and $w=u^*vu$. By direct computation, we get that $w$ is a
unitary element in $\mathcal R_2$
  and
orthogonal to $\mathcal A$ in $\mathcal R_2$ with respect to $\tau$.
\qquad Q.E.D. \vspace{1.3cm}

\section {Some Technical Lemmas}

Let $\{n_i\}_{i=1}^N$ be a sequence of integers where each $n_i\ge
2$ and $N$ is finite or $\infty$. Let $F(n_i)$ be the free group
with the standard generators $\{g_{i,j}\}_{1\le j \le n_i}$.
 Let $$\begin{aligned} G&=\times_{i=1}^N F(n_i), \text { the direct product of groups $F(n_1), \ldots, F(n_N)$};\\
% H&= \text {subgroup of $G$ generated by }\{ g_{i,1}\}_{i=1}^N\\
 H_i &=\text {subgroup of $F(n_i)$ generated by }  g_{i,1}, \qquad \text { for }1\le i\le N\\
 H&= H_1\times H_2\times \cdots \times H_N\\
 G_i& = \left (\times_{k=1}^{i-1}F(n_k)\right ) \times H_i \times \left (\times_{k=i+1}^N
 F({n_k})\right),  \qquad \text { for }1\le i\le N.
\end{aligned}$$
Here, we identify $F(n_i)$ with its canonical image in $G$.

Let $\M_0$ be a finite von Neumann algebra with the tracial state
$\tau_0$. Let $G$ act on $\M_0$ by $\tau_0$-preserving
automorphisms. Denote by $\M=\M_0\times G$ the corresponding crossed
product von Neumann algebra.   $\M_0$ is identified with its
canonical image in $\M$ and denote by $\lambda(g)$, $g \in G$, the
unitary elements in $M$ canonically implementing the action of $G$
on $\M_0$, and by $\tau$  the tracial state on $\M$ that extends
$\tau_0$ of $\M_0$.

\vspace{0.3cm}

Note every element $x$ in $L^2(\M,\tau)$ can be uniquely decomposed
as $x=\sum_{g \in G} a_g\lambda(g)$, with $a_g\in \M_0$. The set
$\{g\in G | a_g \ne 0\}$ is called the support of $x$.

\vspace{0.3cm} If $x=\sum_{g \in G} a_g\lambda(g)$ is in
$L^2(\M,\tau)$ and $F\subset G$ is a nonempty subset of $G$, we
denote by $x_F$ the element $\sum_{g \in F} a_g\lambda(g)\in
L^2(\M,\tau)$ and $||x||_F= ||x_F||_2.$  For any subsets $F$ and
$\tilde F$ of $G$, denote by $F^{-1}$ the set $  \{g^{-1}|
   g\in F\}$ and by $F\tilde F$ the set $\{gh|g\in F, \ h\in \tilde F\}.$

\vspace{0.3cm} For the subgroup $H_0\subset G$, we denote by
$L(H_0)$ the von Neumann subalgebra of $\M$ generated by
$\lambda(g)$ with $g\in H_0$ and by $\M_{H_0}$ the von Neumann
subalgebra of $\M$ generated by $\M_0$ and $L( H_0)$.

\vspace {0.3cm} The following lemma is essentially from Lemma 2.1 in
[9].
\begin{lemma}
Let $\omega$ be a free ultrafilter on $\Bbb N$, and $H$ be a
subgroup of $G$. Suppose $x= (x_n)_n$ is an element in $\M^\omega$
and $y$ is an element in $\M$.  Suppose, for every $\epsilon>0$,
there are subsets $S_0,S_1,  S$ of $G$, depending on $\epsilon$,
satisfying
\begin{enumerate}[(i)]
  \item $S= G\setminus (S_0\cup H)$;
  \item $||y-y_{S_1}||_2 \le  \epsilon,$ and $||y_{S_1}||\le ||y||$;
  \item There exists some positive integer $n_1$ such that
     $  ||(x_n)||_{S_0}\le \epsilon, \quad \forall n\ge n_1;$
   \item  $(SS_1)\cap (S_1S)=\phi ;$  $(S_1H)\cap(S_1S)=\phi;$
  $(HS_1)\cap(S_1S)= \phi;$
   \item $(HS_1)\cap (S S_1)=\phi;$ $(S_1H)\cap (SS_1)= \phi$.
\end{enumerate}
Then, $$ ||yx-xy||_2^2 \ge || y(x-E_{\M_H^{\omega}}(x)) ||_2^2 +
||(x-E_{\M_H^{\omega}}(x))y||_2^2.
$$
\end{lemma}
\noindent {\bf Proof: }  Note that the support of $y_{S_1}$,
$(x_n)_{H}$ or $(x_n)_{S}$ is on $S_1$, $H$ and $S$, respectively.
By (iv) and (v), it is easy to check that $y_{S_1} [(x_n)_{S}]$ ,
$[(x_n)_{S}]y_{S_1}$ and $y_{S_1}[(x_n)_{H}]  - [(x_n)_{H}]
 y_{S_1}$   are mutually orthogonal vectors in $L^2(\M,\tau)$.

\vspace{0.3cm}

Thus, if $\mathcal H_{\omega}$ denote the untraproduct Hilbert space
obtained as the quotient of $(\xi_n)_n \subset L^2(\mathcal M,\tau)|
\ \sup ||\xi_n||_2 < \infty\}$ by the subspace $\{(\eta_n)_n \subset
L^2(\mathcal M,\tau)| \ \lim_{n\rightarrow \omega} ||\xi_n||_2
=0\}$, endowed with the norm $||(\xi_n)_n||_2 =\lim_{n\rightarrow
\omega} ||\xi_n||_2$, then $x'= (y_{S_1}[(x_n)_{S}])_n $, $x''=
([(x_n)_{S}]y_{S_1})_n$, $x''' = (y_{S_1}[(x_n)_{H}]  - [(x_n)_{H}]
 y_{S_1})_n$ are mutually orthogonal elements in $\mathcal
H_{\omega}$. Moreover $L^2(\mathcal M^{\omega}, \tau)$ is naturally
embedded in $\mathcal H_{\omega}$. Note that $E_{\M_H^\omega}(x_n)
=(x_n)_H$ for $n\ge 1$. By (i), (ii), (iii) we have:
$$
\begin{aligned}
  ||y(x-E_{\M_H^{\omega}}(x)) -x'||_2 &\le sup_{n\ge n_1} || y(x_n-E_{\M_H^{\omega}}(x_n)) -  y_{S_1} [(x_n)_{S}]||_2\\
    &\le sup_{n\ge n_1} || (y- y_{S_1})(x_n-E_{\M_H^{\omega}}(x_n))||_2 \\
     &\qquad \qquad \qquad +\sup_{n\ge n_1} ||  y_{S_1}[x_n-E_{\M_H^{\omega}}(x_n)- (x_n)_{S}]||_2\\
    &\le sup_{n\ge n_1} || (y- y_{S_1})(x_n-E_{\M_H^{\omega}}(x_n))||_2\\
    & \qquad \qquad \qquad +||y||(\sup_{n\ge n_1} || x_n- (x_n)_H- (x_n)_{S}||_2)\\
    &\le sup_{n\ge n_1} || (y- y_{S_1})(x_n-E_{\M_H^{\omega}}(x_n))||_2\\
    & \qquad \qquad \qquad +||y||(\sup_{n\ge n_1} ||  (x_n)_{S_0}||_2)\\
  &\le \epsilon \sup (||x_n|| + ||y||)\\
   \end{aligned}$$$$\begin{aligned}
  ||(x-E_{\M_H^{\omega}}(x))y -x''||_2 &\le sup_{n\ge n_1} ||  (x_n-E_{\M_H^{\omega}}(x_n))y -  y_{S_1} [(x_n)_{S}]||_2\\
    &\le sup_{n\ge n_1} || (x_n-E_{\M_H^{\omega}}(x_n))(y- y_{S_1})||_2\\
     & \qquad \qquad \qquad+(\sup_{n\ge n_1} ||     (x_n)_{S_0} ||_2)||y||\\
  &\le \epsilon \sup (||x_n|| + ||y||)\\
 ||y E_{\M_H^{\omega}}(x)  -  E_{\M_H^{\omega}}(x)y- x''' ||_2 &\le sup_{n\ge n_1} ||
 (y- y_{S_1})E_{\M_H^{\omega}}(x_n)
 ||_2\\
  & \qquad \qquad \qquad+ sup_{n\ge n_1} || E_{\M_H^{\omega}}(x_n)(y- y_{S_1})||_2\\
  &\le2 \epsilon \sup ||x_n||
\end{aligned}
$$
This shows that the vectors $y(x-E_{\M_H^{\omega}}(x))$,
$(x-E_{\M_H^{\omega}}(x))y$, $y E_{\M_H^{\omega}}(x) -
E_{\M_H^{\omega}}(x)y$ can be approximated arbitrarily well in
$\mathcal H_{\omega}$ by some mutually orthogonal vectors and hence
they are mutually orthogonal in $L^2(\mathcal M^{\omega}, \tau)$.
Since there sum is equal to $yx- xy$ we get
$$\begin{aligned}
  ||yx-xy||_2^2& = ||y(x-E_{\M_H^{\omega}}(x))||_2^2 + ||(x-E_{\M_H^{\omega}}(x))y||_2^2
   +||yE_{\M_H^{\omega}}(x)  -  E_{\M_H^{\omega}}(x)y||_2^2
  \\
  &\ge ||y(x-E_{\M_H^{\omega}}(x))||_2^2 + ||(x-E_{\M_H^{\omega}}(x))y||_2^2
  \end{aligned}
$$
Q.E.D. \vspace{1.5cm}

\vspace{0.3cm} Recall that, for the subgroup $H_0\subset G$, we
denote by $L(H_0)$ the von Neumann subalgebra of $\M$ generated by
$\lambda(g)$ with $g\in H_0$ and by $\M_{H_0}$ the von Neumann
subalgebra of $\M$ generated by $\M_0$ and $L( H_0)$. Following the
preceding notations, we let
$$
\A=\M_{H},\quad \A_i=\M_{H_i},\quad \N_i= \M_{G_i}, \qquad \text {
for } 1\le i\le N.
$$

 \pb We have the following lemma, which is the extension of Lemma
 2.1 in [9].
\begin{lemma} Suppose that $N$ is a finite integer. Let $\omega$ be a free ultrafilter on $\Bbb N$.
 Suppose $x$ is an element in $\M^{\omega}$ $(=(\M_0\times G)^\omega)$ that commutes
 $\mathcal A$ and $$E_{\N_1^{\omega}}(x)= \cdots =E_{\N_N^\omega}(x) =E_{\A^\omega}(x) .$$ Then for any
$y \in \M$ with $ E_{\N_1}(y)= \cdots =E_{\N_N}(y) =E_{\A }(y) = 0$,
we have
$$
||yx-xy||_2^2 \ge || y(x-E_{\mathcal A^{\omega}}(x)) ||_2^2 +
||(x-E_{\mathcal A^{\omega}}(x))y||_2^2.
$$
\end{lemma}

\noindent {\bf Proof: } Let $(x_n)_n$ be a sequence of elements in
$\mathcal M$ representing $x \in \mathcal M^{\omega}$. We might
assume that
$$ E_{\N_1 }(x_n)= \cdots =E_{\N_N}(x_n)=E_{\A}(x_n), \ \forall n \in \Bbb
N. $$

\pb  Let $\mathcal F = span \{ a_g\lambda (g)| a_g\in \M_0, g \in
G\}$, a weakly dense *-subalgebra in $\mathcal M$. For every
$\epsilon
>0$, by Kaplansky density theorem there exists $y^0
 \in \mathcal F$ such that $$||y-y^0||_2 < \epsilon, \
  \   ||y^0||\le ||y||,  \   E_{\N_1 }(y^0)= \cdots =E_{\N_N}(y^0)= E_{\mathcal A}(y^0)=0.
$$
Let  $S_1$ be the support of $y^0$. Since $E_{\N_1}(y^0)=\cdots =
E_{\N_N}(y^0)=E_{\mathcal A}(y^0)=0$, we have that $S_1\cap (\cup_i
G_i\cup H)=\phi$.

  Note that every element (or word) $w$ in $G=F(n_1)\times\cdots
\times  F(n_N)$ can be uniquely written as $$w=
(g_{1,1})^{m_1}\cdots (g_{N,1})^{m_N} w_1\cdots
w_N(g_{1,1})^{n_1}\cdots (g_{N,1})^{n_N},$$ where
$(g_{i,1})^{m_i}w_i(g_{i,1})^{n_i}$ is a reduced word in $F(n_i)$
for $1\le i \le N$.

Let $N_0-1$ be the maximal length of the words $g$ in the finite set
$S_1$, the supports of $y^0$. For every $i$, denote by $S_i^0 =
\{g=(g_{1,1})^{m_1}\cdots (g_{N,1})^{m_N} w_1\cdots
w_N(g_{1,1})^{n_1}\cdots (g_{N,1})^{n_N}\in G | \   w_i $ starts
with   a nonzero power of some $g_{i,j}$  for some $ j\ge 2$; and
$0\le |m_i|\le 2N_0-1\}$,  and
$$\begin{aligned} S_0 &= (\cup_iS_i^0\cup (\cup_iS_i^0)^{-1} \cup G_1\cup
G_2\cup \cdots\cup G_N)\setminus H,\\ S&= G \setminus (S_0\cup H).
\end{aligned}$$

Our first goal is to show that $||(x_n)_{S_0}||_2$ is small for $n$
large. Note that $(x_n)_{G_i}=E_{\N_i}(x_n)= E_{\A}(x_n)=(x_n)_H$,
for $1 \le i \le N$. It follows that $||(x_n)_{S_0}||_2 \le
\sum_i||(x_n)_{S_i^0}||_2 +\sum_i||(x_n)_{(S_i^0)^{-1}}||_2 $. It
will be sufficient to control the norms in the right side. Let $N_1$
be an integer multiple of $4N_0$ such that $N_1\ge 32N_0N^3||x||^2
{\epsilon^{-2}} $. By hypothesis, there exists $n_1 = n_1(\epsilon,
N_1)$ such that if $n\ge n_1$, then
$$
||\lambda(g_{i,1})^{k_i}  x_n \lambda(g_{i,1})^{-k_i} -x_n||_2 <
(2N)^{-2} \epsilon;
$$ for all $1\le i \le N$ , $|k_1|, \ldots, |k_N| \le N_1$. So if $1\le i\le N, 0<4N_0|k_i|  \le N_1$ and $n\ge
n_1$, then we have
$$
 \begin{aligned}
     ||\lambda&(g_{i,1})^{4N_0k_i}  (x_n)_{S_i^0}\lambda(g_{i,1})^{-4N_0k_i}
     -
     (x_n)_{  (g_{i,1})^{4N_0k_i}   S_i^0 (g_{i,1})^{-4N_0k_i}}||_2\\
     &= ||  (\lambda (g_{i,1})^{4N_0k_i}  x_n\lambda(g_{i,1})^{-4N_0k_i}    -x_n)_{  (g_{i,1})^{4N_0k_i}  S_i^0
      (g_{i,1})^{-4N_0k_i}}||_2\\
     &\le ||  \lambda (g_{i,1})^{4N_0k_i} x_n
    \lambda(g_{i,1})^{-4N_0k_i}
     -x_n ||_2 \\&< (2N)^{-2} \epsilon;
 \end{aligned}
$$
Using the parallelogram identity in the Hilbert space $L^2(\mathcal
M,\tau)$, we get the inequalities
$$
\begin{aligned}
  ||(x_n)_{S_i^0}||_2^2 &=
  ||\lambda(g_{i,1})^{4N_0k_i}  (x_n)_{S_i^0}\lambda(g_{i,1})^{-4N_0k_i}||_2^2\\
  &\le 2  ||\lambda(g_{i,1})^{4N_0k_i}  (x_n)_{S_i^0}\lambda(g_{i,1})^{-4N_0k_i} -
     (x_n)_{ (g_{i,1})^{4N_0k_i}  S_i^0
      (g_{i,1})^{-4N_0k_i} }||_2^2 \\
     & \qquad \qquad \quad + 2|| (x_n)_{ (g_{i,1})^{4N_0k_i}  S_i^0
      (g_{i,1})^{-4N_0k_i}} ||_2^2\\ & \le
     (2N)^{-3}\epsilon^2
     +2 || (x_n)_{ (g_{i,1})^{4N_0k_i}  S_i^0
      (g_{i,1})^{-4N_0k_i}} ||_2^2.
\end{aligned}
$$
%Similarly,
%$$
%\begin{aligned}
%  ||(x_n)_{S_0^0}||_2^2   \le
%     2^{-3}\epsilon
%     +2 || (x_n)_{\lambda(g_1)^{4N_0k_1}\lambda(h_1)^{ 4N_0k_2}S_0^0\lambda(h_1)^{- 4N_0k_2}\lambda(g_1)^{-4N_0k_1}} ||_2^2.
%\end{aligned}
%$$
Now we use the fact that $$\{ \lambda (g_{i,1})^{4N_0k_i} S_i^0
     \lambda (g_{i,1})^{-4N_0k_i}
\}_{k_i\in \Bbb Z}$$ are disjoint subsets of $G$, so that summing up
the above inequalities for all $k_i$, $0< 4 N_0|k_i| \le N_1$, we
have
$$
  \left ( \frac {N_1}{2N_0}  \right )  ||(x_n)_{S_i^0}||_2^2 < \left ( \frac {N_1}{2N_0}  \right
  )   (2N)^{-3}\epsilon^2 + 2 ||x_n||_2^2
$$
so that
$$
   ||(x_n)_{S_i^0}||_2^2 < {(2N)}^{-3}\epsilon^2 + 2||x||^2  \left ( \frac {N_1}{2N_0}  \right
   )^{-1} \le {(2N)}^{-2} \epsilon^2.
$$
Similarly, we get $||(x_n)_{(S_i^0)^{-1}}||_2 < {(2N)}^{-1}\epsilon$
and thus $||(x_n)_{S_0}||_2 <  \epsilon$ for all $n \ge n_1$.

\vspace{0.3cm}

Denote by $l(w)$ the length of the reduced word $w$ in $G$. Since
every element $w$ in $G$ can be uniquely expressed as,
$$w=(g_{1,1})^{m_1}\cdots (g_{N,1})^{m_N} w_1\cdots
w_N(g_{1,1})^{n_1}\cdots (g_{N,1})^{n_N}\in G,$$ where
$g_{i,1}^{m_i}w_ig_{i,1}^{n_i}$ is a reduced word in $L(F(n_i))$ and
$w_i$ does not start, also not end, with any power of $g_{i,1}$.
Then we can define the following functions as
$$
Sl_i(w)=m_i,\quad El_i(w)= n_i,\quad l_i(w)= l(w_i),\quad \tilde
l(w) =l(w_1\cdots w_N).
$$

Note that every reduced word $w$ in $S$ contains, for all $1\le i\le
N$, a nonzero power of some $g_{i,j}$   with $j\ge 2$; and begins,
also ends, with the power of   $g_{i,1}$   greater than in $2N_0-1$.
i.e.
$$
\min_il_i(w)>0;\quad \min_{i}Sl_i(w)\ge 2N_0-1;\quad
\min_{i}El_i(w)\ge 2N_0-1.
$$

Let $g_1$ be any element in $S_1$, the support of $y^0$. Since
$S_1\cap (\cup_{1\le i \le N}G_i\cup H)=\phi$, it follows that, for
every $i$, $g_1 $
  contains a nonzero power of some
$g_{i,j}$ $(j\ge 2)$. (i.e. $\min_i l_i(g_1)>0$). Since any element
in $S$ begins, also ends, with the powers of $g_{i,1}$ greater than
in absolute value that twice the length of $g_1$, a quick
computation shows that
$$ \begin{aligned} &\min_{i}Sl_i(w)\ge 2N_0-1;\quad  \min_iEl_i(w)\le N_0, \quad
\forall w \in SS_1\\
&\min_{i}Sl_i(w)\le N_0;\quad  \min_iEl_i(w) \ge 2N_0-1, \quad
\forall w \in  S_1S\end{aligned}
$$ Hence,
$ (SS_1)\cap(S_1S)=\phi .$

Let $g_2$ be another element in $S_1$ (so, not in $\cup_{1\le i \le
N}G_i\cup H$). We claim that $g_2S\cap g_1 H = \phi$. Since $g_2$ is
in $S_1$ but not in $\cup_{1\le i \le N}G_i\cup H$, we know that
$\min_il_i(g_2)>0$  and $l(g_2)\le N_0$. Combining with the facts
that $\min_{i}Sl_i(g)\ge 2N_0-1$  and $ \min_i l_i(g)>0$ for every
$g\in S$, we obtain that $l_i(g_2g)\ge N_0+1$.
 Therefore, $\tilde l(g_2g)\ge l_i(g_2g)\ge  N_0+1$. On the other hand, it is
easy to see that $\tilde l(g_1h)\le N_0$ for all $g_1\in S_1$ and
$h\in H$.   Thus $g_2S\cap g_1 H = \phi$ whence
$(S_1H)\cap(S_1S)=\phi$. Similarly, we also have
  $(HS_1)\cap(S_1S)= \phi$;
    $(HS_1)\cap (S S_1)=\phi;$ and $(S_1H)\cap(SS_1)= \phi$.

As a summary, for such $x,y, H, S_0,S_1,S$, we have that (i) $S=
G\setminus (S_0\cup H)$;
  (ii) $||y-y_{S_1}||_2 \le  \epsilon,\ ||y_{S_1}||\le ||y||;$
  (iii) There exists some positive integer $n_1$ such that
     $  ||(x_n)||_{S_0}\le \epsilon, \quad \forall n\ge n_1;$
  (iv) $(SS_1)\cap (S_1S)=\phi $;  $(S_1H)\cap(S_1S)=\phi$;
  $(HS_1)\cap(S_1S)= \phi$;
   (v) $(HS_1)\cap (S S_1)=\phi;$ $(S_1H)\cap(SS_1)= \phi$.

   Applying Lemma 3.1, we get
   $$
||yx-xy||_2^2 \ge || y(x-E_{\mathcal A^{\omega}}(x)) ||_2^2 +
||(x-E_{\mathcal A^{\omega}}(x))y||_2^2.\qquad \qquad\text {\qquad
Q.E.D}
$$

\pb Note when $N=1$, $G_1=H_1=H$, as an application of the preceding
lemma, we have the following corollary.

\begin{Col}
 Let $\M_0$ be $L(\Bbb Z)$ or $\Bbb CI$, and the group $G$ $(=F(n ))$ act trivially on $\M_0$. Then
\begin{enumerate}
\item  $\M=\M_0\times G=
 L(G) \otimes\M_0$, and $\A=L(H) \otimes \M_0$.
\item  Let $\omega$ be a free ultrafilter on $\Bbb N$.
 Suppose $x$ is an element in $\M^{\omega}$ $(=(\M_0\times G)^\omega)$ that commutes
 $\mathcal A$.   Then for any
$y \in \M$ with $E_{\A }(y) = 0$, we have
$$
||yx-xy||_2^2 \ge || y(x-E_{\mathcal A^{\omega}}(x)) ||_2^2 +
||(x-E_{\mathcal A^{\omega}}(x))y||_2^2.
$$\end{enumerate}
\end{Col}

\pb Generally, when $N$ is an arbitrary finite integer, we have the
following corollary.
\begin{Col}
 Let $\M_0$ be $L(\Bbb Z)$ or $\Bbb CI$, and the group $G$ $(=F(n_1)\times \cdots
 \times F(n_N))$ act trivially on $\M_0$. Then
\begin{enumerate}
\item  $\M=\M_0\times G=
 L(G) \otimes\M_0$, and $\A=L(H) \otimes \M_0$.
\item  Suppose that $N$ is a finite integer. Let $\omega$ be a free ultrafilter on $\Bbb N$.
 Suppose $x$ is an element in $\M^{\omega}$ $(=(\M_0\times G)^\omega)$ that commutes
 $\mathcal A$ and $$E_{(L(G_1)\otimes \M_0)^{\omega}}(x)= \cdots =E_{(L(G_N)\otimes \M_0)^\omega}(x) =E_{\A^\omega}(x) .$$ Then for any
$y \in \M$ with $ E_{L(G_1)\otimes \M_0}(y)= \cdots
=E_{L(G_N)\otimes \M_0}(y) =E_{\A }(y) = 0$, we have
$$
||yx-xy||_2^2 \ge || y(x-E_{\mathcal A^{\omega}}(x)) ||_2^2 +
||(x-E_{\mathcal A^{\omega}}(x))y||_2^2.
$$\end{enumerate}
\end{Col}

\vspace{2cm}

\section{  Abelian Maximal Injective Subalgebras of Tensor Product of  Free Group Factors  }

Let $\{n_i\}_{i=1}^m$ be a sequence of integers where $n_i\ge 2$ for
all $1\le i\le m$ and $m$ is finite or infinite. Let $F(n_i)$ be the
free group with the standard generators $\{g_{i,j}\}_{j=1}^{n_i}$
for all $1\le i\le m$.
 Let $$\begin{aligned} G(m)&=\times_{i=1}^m F(n_i), \text { the direct product of groups $F(n_1), \ldots, F(n_m)$};\\
% H(m)&= \text {subgroup of $G$ generated by }\{ g_{i,1}\}_{i=1}^m;\\
 H_i &=\text {subgroup of $F(n_i)$ generated by }  g_{i,1}, \qquad \text { for }1\le i\le m;\\
 H(m)&= H_1\times H_2\times \cdots \times H_m\\
 G_i& = \left (\times_{k=1}^{i-1}F(n_k)\right ) \times H_i \times \left (\times_{k=i+1}^m
 F_{n_k}\right),  \qquad \text { for }1\le i\le m;\\
 J_i &=\left (\times_{k=1}^{i-1}F(n_k)\right )   \times \left (\times_{k=i}^m
 H_{k}\right),  \qquad \text { for }1\le i\le m.\\
\end{aligned}$$
 Here, $F(n_i)$ is identified with its canonical image in $G(m)$.

\pb Let $\lambda$ be the left regular representation of $G(m)$ and
$\M(m)=L(G(m))$ the group von Neumann algebra associated with
$G(m)$. Denote by $\A(m)$ (or $\A_i$, $\N_i$, $\B_i$) the von
Neumann subalgebra $L(H(m))$, (or $L(H_i)$, $L(G_i)$, $L(J_i)$
respectively) of $\M$, for all $1\le i\le m$. \vspace{0.3cm}

  From the construction of $\A(m)$, it is easy to see
the following lemma.

\begin{lemma}
  $\mathcal A(m)$ is a maximal abelian   von Neumann subalgebra of
  $\mathcal M(m)$.
\end{lemma}

\vspace{0.3cm} Moveover, we have that
\begin{lemma}
 $\mathcal A(m)$ is a singular maximal abelian   von Neumann subalgebra  of $
 \mathcal M(m)$.
\end{lemma}
\noindent {\bf Proof: }   For any group element $g $   in $G(m)
\setminus (\cup_i G_i)$, we have $\lambda (g) \mathcal A(m) \lambda
(g)^* $ and $\mathcal A(m)$ are mutually orthogonal in $\M(m)$. By
Lemma 2.5 in [7], it follows that $\lambda (g)$ is orthogonal to
$\N(\mathcal A(m))$, where $\N(\mathcal A(m))$ is the von Neumann
subalgebra of $\M(m)$ generated by the normalizer of $\mathcal A(m)$
in $\M(m)$.
 But the Hilbert space generated in $L^2(G(m))$ by $\lambda(g)$ with
 $g\in G(m)\setminus (\cup_i G_i)$, coincides with the orthogonal of
 $L^2(\cup  G_i)$ in $L^2(G(m))$.
Therefore $\N(\mathcal A(m))\subset  L^2(\cup_i G_i)$. Moreover, for
all $i\ne j$, if the unitary element $\lambda(g)$ is an element such
that $g$ is in $G_i$ and not contained in $G_j$, then $\lambda
(g)\mathcal A_j \lambda (g)^*$ and $\mathcal A(m)$ are mutually
orthogonal. Hence, again by Lemma 2.5 in [7], such $\lambda (g)$ is
orthogonal to $\N(\mathcal A(m))$. Therefore, for all $g$ not
contained in $\cap_iG_i$, $\lambda(g)$ is orthogonal to $\N(\A(m))$.
It follows that $ \N(\mathcal A(m))=\A(m)$.
 Combining with the preceding lemma, $\A(m)$ is a singular maximal abelian   von
Neumann subalgebra of $\M(m)$. \qquad Q.E.D.

\vspace{0.3cm}

 \begin{lemma}Let $\B\cong L(\Bbb Z)$. Then we also have that $\A(m) \otimes \B$ is a singular maximal
  abelian   von Neumann subalgebra of $\M(m)\otimes \B$.
 \end{lemma}
\noindent{Proof: } The proof of the lemma is almost identical with
the one of Lemma 4.2. So we skipped it here. Q.E.D.

\vspace{0.3cm} The following lemma is Corollary 3.3 in [9]. For the
reader's convenience, we present a proof here.

\begin{lemma}
$\A(1)$ is a maximal injective abelian von Neumann subalgebra in
$\M(1)$.
\end{lemma}
\noindent {\bf Proof: } Assume $\mathcal R$ is a maximal injective
von Neumann subalgebra of $\M(1) $ and $\mathcal A(1)  \subset
\mathcal R\subset \mathcal M(1)  $.
%Suppose $\mathcal R= \mathcal R_1 \oplus
%\mathcal R_2$, where $R_1$ is a type I von Neumann subalgebra of
%$\mathcal M(k)  $ and $\mathcal R_2$ is a type II$_1$ injective von
%Neumann subalgebra of $\mathcal M(k)  $. If $\mathcal R_2 \ne 0$,
%from Corollary 2.1 we can find some $x$ in $\mathcal R'\cap \mathcal
%R_2^{\omega}$ but not contained in $\mathcal A(k)  ^{\omega}$. By
%Lemma 2.4 we can find a unitary $w$ in $\mathcal R_2$ such that $w$
%is orthogonal to $\mathcal A(k)  $ in $\R$.
Let $\M_0=\Bbb C$ and the group $G(1)$ act trivially on $\M_0$.

Decompose $\mathcal R = \mathcal R_1\oplus \mathcal R_2$ where
$\mathcal R_1$ is a type I von Neumann subalgebra and $\mathcal R_2$
is a type II$_1$ injective von Neumann subalgebra. If $\mathcal R_2
\ne 0$, from Corollary 2.1 we can find some $x$ in $\mathcal R'\cap
\mathcal R_2^{\omega}$ but not contained in $\mathcal
A(1)^{\omega}$. By Lemma 2.4 we can find a unitary $w$ in $\mathcal
R_2$ such that $w$ is orthogonal to $\mathcal A(1)$ in $\R$, whence
$ E_{\A(1)}(w)=0$. By Corollary 3.1 and the fact that $x\in
\R_2^\omega$, we have that $0 =||xw-wx||_2\ge || (x-E_{\mathcal
A(1)^{\omega}}(x))w||_2
>0 $, which is a contradiction. Therefore $\mathcal R_2=0 $ and
$\mathcal R=\mathcal R_1$. From Lemma 2.3 and Lemma 4.2, it follows
that $\mathcal A(1)  =\mathcal R$. \quad Q.E.D

\begin{lemma} Let $\B\cong L(\Bbb Z)$. $\A(1)\otimes \B$ is a maximal injective
subalgebra of $\M(1)\otimes \B$ $(=L(F(n_1))\otimes \B)$.
\end{lemma}

 \noindent {\bf Proof: }
 Suppose $\mathcal R$ is a
maximal injective von Neumann subalgebra and $\mathcal A(1)\otimes
\B \subset \mathcal R\subset \M(1)\otimes \B$. Suppose $\mathcal R=
\mathcal R_1 \oplus \mathcal R_2$, where $\mathcal R_1$ is a type I
von Neumann subalgebra of $\M(1)\otimes \B$ and $\mathcal R_2$ is a
type II$_1$ injective von Neumann subalgebra of $\M(1)\otimes \B$.
If $\mathcal R_2 \ne 0$, from Corollary 2.1 we can find some $x$ in
$\mathcal R'\cap \mathcal R_2^{\omega}$ but not contained in
$(\mathcal A(1)\otimes \B)^{\omega}$. By Lemma 2.4 we can find a
unitary $w$ in $\mathcal R_2$ such that $w$ is orthogonal to
$\mathcal A(1)\otimes \B$ in $\R$. Let $\M_0=\B$ and the group
$G(1)$ act trivially on $\B$. By Corollary 3.1 and the fact that $x
\in \R_2^{\omega}$, we have that $0 =||xw-wx||_2 \ge
||(x-E_{(\mathcal A(1)\otimes \B)^{\omega}}(x))w||_2 =
||x-E_{(\mathcal A(1)\otimes \B)^{\omega}}(x) ||_2>0 $, which is a
contradiction. Therefore $\mathcal R_2=0 $ and $\mathcal R=\mathcal
R_1$. From Lemma 2.3 and Lemma 4.2, it follows that $\mathcal
A(1)\otimes \B=\mathcal R$. \qquad Q.E.D.

\begin{lemma}Let $\B\cong L(\Bbb
Z)$.  Suppose that, when $m< k$, $\A(m)\otimes \B$ is a maximal
injective von Neumann subalgebra  of $\M(m)\otimes \B$. Then
$\A(k)\otimes \B$ is also a maximal injective von Neumann subalgebra
of $\M(k)\otimes \B$.

\end{lemma}
\noindent {Proof: } Assume $\mathcal R$ is a maximal injective von
Neumann subalgebra of $\M(k)\otimes \B$ and $\mathcal A(k)\otimes \B
\subset \mathcal R\subset \mathcal M(k)\otimes \B$. %Suppose
%$\mathcal R= \mathcal R_1 \oplus \mathcal R_2$, where $R_1$ is a
%type I von Neumann subalgebra of $\mathcal M(k)\otimes \B$ and
%$\mathcal R_2$ is a type II$_1$ injective von Neumann subalgebra of
%$\mathcal M(k)\otimes \B$. If $\mathcal R_2 \ne 0$, from Corollary
%2.1 we can find some $x$ in $\mathcal R'\cap \mathcal R_2^{\omega}$
%but not contained in $(\mathcal A(k)\otimes \B)^{\omega}$. By Lemma
%2.4 we can find a unitary $w$ in $\mathcal R_2$ such that $w$ is
%orthogonal to $\mathcal A(k)\otimes \B$ in $\R$.
Let $\M_0=\B$ and
the   group $G(k)$ act trivially on $\M_0$.

  Since, for all $1\le i\le k$, $(L(G_i)\otimes \B)'\cap (\mathcal M(k)\otimes \B)\subset (L(G_i)\otimes \B)$,
there is a unique trace-preserving  condition expectation
$E_{L(G_i)\otimes \B}$ from $\mathcal M(k)\otimes \B$ onto $
L(G_i)\otimes \B$. Actually $
   E_{L(G_i)\otimes \B}(x)$  is defined as $$E_{L(G_i)\otimes \B}(x)=\lim_{n \rightarrow
\infty}
   \frac 1 {2n} \sum_{l=-n}^n \lambda(g_{i,1})^l x \lambda(g_{i,1})^{-l}.
$$ Moreover, if $x$ is expressed as $\sum_{g \in G}a_g\lambda (g)$,
then   $E_{ L(G_i)\otimes \B}(x)$ $= \sum_{g \in G_i}a_g\lambda
(g),$ where $a_g \in \B$.

If $E_{L(G_i)\otimes \B}(\mathcal R) \supsetneqq \mathcal
A(k)\otimes \B$, there exists some $x$ in $\mathcal R$ such that
$E_{
 L(G_i)\otimes \B} (x) $ is not contained in $\mathcal A(k)\otimes \B$.
From the fact that $ E_{L(G_i)\otimes \B}(x) =\lim_{n \rightarrow
\infty}
   \frac 1 {2n} \sum_{l=-n}^n \lambda(g_{i,1})^l x \lambda(g_{i,1})^{-l}$ and $\lambda(g_{i,1})$ belongs
   to $\A_i\subset\A(k)$, we get that $ E_{L(G_i)\otimes \B}(x)$ is also contained in $\mathcal R$.
   Denote  $
E_{L(G_i)\otimes \B}(x)$  by $y$. Therefore $y$ is in $\mathcal
R\cap (L(G_i)\otimes \B)$ but not contained in $\mathcal A(k)\otimes
\B$. Let $\mathcal S$ be the von Neumann subalgebra generated by $y$
and $\mathcal A(k)\otimes \B$ in $\mathcal R\cap (L(G_i)\otimes
\B)$. Since $\mathcal R$ is injective, $\mathcal S$ is also
injective and contained in $L(G_i)\otimes \B$. Note that
$$\begin{aligned} \A(k)\otimes \B &=\left (\otimes_{j=1}^{i-1}\A_j\right ) \otimes \left (\otimes_{j=i+1}^k\A_j\right ) \otimes
\left (\A_i\otimes \B\right )\\& \subset \mathcal S\subset \left
(\otimes_{j=1}^{i-1}L(F(n_j))\right ) \otimes \left
(\otimes_{j=i+1}^kL(F(n_j))\right ) \otimes \left (\A_i\otimes
\B\right )\end{aligned}$$
 By induction
hypothesis, we know $ \mathcal S= \mathcal A(k)$, contradicting with
the assumption that $0\ne y\in \mathcal S\setminus \A(k)$. Hence we
obtain $E_{ L(G_i)\otimes \B }(\mathcal R) = \mathcal A(k)\otimes
\B$ for all $1\le i\le k$. It follows that we have $E_{L(G_1)\otimes
\B}(\mathcal R) =\cdots=E_{L(G_k)\otimes \B}(\R)= \mathcal
A(k)\otimes \B$.

Therefore, we can assume that $E_{L(G_1)\otimes \B}(\mathcal R)
=\cdots=E_{L(G_k)\otimes \B}(\R)=  {\mathcal A(k)\otimes \B}$. Again
decompose $\mathcal R = \mathcal R_1\oplus \mathcal R_2$ where
$\mathcal R_1$ is a type I von Neumann subalgebra and $\mathcal R_2$
is a type II$_1$ injective von Neumann subalgebra. If $\mathcal R_2
\ne 0$, from Corollary 2.1 we can find some $x$ in $\mathcal R'\cap
\mathcal R_2^{\omega}$ but not contained in $(\mathcal A(k)\otimes
B)^{\omega}$.   By Lemma 2.4 we can find a unitary $w$ in $\mathcal
R_2$ such that $w$ is orthogonal to $\mathcal A(k)\otimes \B$ in
$\R$, whence $E_{L(G_1)\otimes \B}(w) =\cdots=E_{L(G_k)\otimes
\B}(w)
 =E_{\A(k)\otimes \B}(w)=0$. By Corollary 3.2 and the fact that $x\in
\R_2^\omega$, we have that $0 =||xw-wx||_2\ge || (x-E_{\mathcal
A(k)^{\omega}}(x))w||_2
>0 $, which is a contradiction. Therefore $\mathcal R_2=0 $ and
$\mathcal R=\mathcal R_1$. From Lemma 2.3 and Lemma 4.2, it follows
that $\mathcal A(k)\otimes \B=\mathcal R$.\quad Q.E.D.

%\vspace{3cm}

%By Lemma 2.1 in [6] and the fact that $x \in \R_2^{\omega}$, we have
%that $0 =||xw-wx||_2 \ge ||(x-E_{(\mathcal A\otimes
%\B)^{\omega}}(x))w||_2 = ||(x-E_{(\mathcal A\otimes
%\B)^{\omega}}(x)) ||_2>0 $, which is a contradiction. Therefore
%$\mathcal R_2=0 $ and $\mathcal R=\mathcal R_1$. From Lemma 2.3 and
%4.3, it follows that $\mathcal A\otimes \B=\mathcal R$. Hence
%$\A\otimes \B$ is maximal injective in $\M(k)\otimes \B$. Q.E.D.

\begin{lemma}
Let $\B\cong L(\Bbb Z)$. Assume that, when $m<k$, $A(m)\otimes \B$
is a maximal injective von Neumann subalgebra in $\M(m)\otimes \B$.
Then $A(k)$ is a maximal injective von Neumann subalgebra of
$\M(k)$.
\end{lemma}

\noindent {Proof: } Assume $\mathcal R$ is a maximal injective von
Neumann subalgebra of $\M(k) $ and $\mathcal A(k)  \subset \mathcal
R\subset \mathcal M(k)  $.
%Suppose $\mathcal R= \mathcal R_1 \oplus
%\mathcal R_2$, where $R_1$ is a type I von Neumann subalgebra of
%$\mathcal M(k)  $ and $\mathcal R_2$ is a type II$_1$ injective von
%Neumann subalgebra of $\mathcal M(k)  $. If $\mathcal R_2 \ne 0$,
%from Corollary 2.1 we can find some $x$ in $\mathcal R'\cap \mathcal
%R_2^{\omega}$ but not contained in $\mathcal A(k)  ^{\omega}$. By
%Lemma 2.4 we can find a unitary $w$ in $\mathcal R_2$ such that $w$
%is orthogonal to $\mathcal A(k)  $ in $\R$.
Let $\M_0=\Bbb C$ and the group $G(k)$ act trivially on $\M_0$.

  Since, for all $1\le i\le k$, $L(G_i)  '\cap \mathcal M(k)  \subset L(G_i)  $,
there is a unique trace-preserving  condition expectation $E_{L(G_i)
}$ from $\mathcal M(k)  $ onto $ L(G_i)  $. Actually $
   E_{L(G_i)  }(x)$  is defined as $$E_{L(G_i)  }(x)=\lim_{n \rightarrow
\infty}
   \frac 1 {2n} \sum_{l=-n}^n \lambda(g_{i,1})^l x \lambda(g_{i,1})^{-l}.
$$ Moreover, if $x$ is expressed as $\sum_{g \in G}a_g\lambda (g)$,
then   $E_{ L(G_i)  }(x)$ $= \sum_{g \in G_i}a_g\lambda (g),$ where
$a_g \in \M_0=\Bbb CI$.

If $E_{L(G_i)  }(\mathcal R) \supsetneqq \mathcal A(k)  $, there
exists some $x$ in $\mathcal R$ such that $E_{
 L(G_i)  } (x) $ is not contained in $\mathcal A(k)  $.
From the fact that $ E_{L(G_i)  }(x) =\lim_{n \rightarrow \infty}
   \frac 1 {2n} \sum_{l=-n}^n \lambda(g_{i,1})^l x \lambda(g_{i,1})^{-l}$ and $\lambda(g_{i,1})$ belongs
   to $\A_i\subset\A(k)$, we get that $ E_{L(G_i)  }(x)$ is also contained in $\mathcal R$.
   Denote  $
E_{L(G_i)  }(x)$  by $y$. Therefore $y$ is in $\mathcal R\cap L(G_i)
$ but not contained in $\mathcal A(k)  $. Let $\mathcal S$ be the
von Neumann subalgebra generated by $y$ and $\mathcal A(k) $ in
$\mathcal R\cap L(G_i)  $. Since $\mathcal R$ is injective,
$\mathcal S$ is also injective and contained in $L(G_i)  $. Note
that
$$\begin{aligned} \A(k)   &=\left (\otimes_{j=1}^{i-1}\A_j\right ) \otimes \left (\otimes_{j=i+1}^k\A_j\right ) \otimes
 \A_i   \\& \subset \mathcal S\subset \left
(\otimes_{j=1}^{i-1}L(F(n_j))\right ) \otimes \left
(\otimes_{j=i+1}^kL(F(n_j))\right ) \otimes  \A_i
\end{aligned}$$
 By induction
hypothesis, we know $ \mathcal S= \mathcal A(k)$, contradicting with
the assumption that $0\ne y\in \mathcal S\setminus \A(k)$. Hence we
obtain $E_{ L(G_i)   }(\mathcal R) = \mathcal A(k)$ for all $1\le
i\le k$. It follows that we have $E_{L(G_1)  }(\mathcal R)
=\cdots=E_{L(G_k)  }(\R)= \mathcal A(k)$.

Therefore, we can assume that $E_{L(G_1)  }(\mathcal R)
=\cdots=E_{L(G_k)  }(\R)= {\mathcal A(k)}$.  Again decompose
$\mathcal R = \mathcal R_1\oplus \mathcal R_2$ where $\mathcal R_1$
is a type I von Neumann subalgebra and $\mathcal R_2$ is a type
II$_1$ injective von Neumann subalgebra. If $\mathcal R_2 \ne 0$,
from Corollary 2.1 we can find some $x$ in $\mathcal R'\cap \mathcal
R_2^{\omega}$ but not contained in $\mathcal A(k)^{\omega}$.   By
Lemma 2.4 we can find a unitary $w$ in $\mathcal R_2$ such that $w$
is orthogonal to $\mathcal A(k)$ in $\R$, whence $E_{L(G_1) }(w)
=\cdots=E_{L(G_k) }(w)
 =E_{\A(k)}(w)=0$. By Corollary 3.2 and the fact that $x\in
\R_2^\omega$, we have that $0 =||xw-wx||_2\ge || (x-E_{\mathcal
A(k)^{\omega}}(x))w||_2
>0 $, which is a contradiction. Therefore $\mathcal R_2=0 $ and
$\mathcal R=\mathcal R_1$. From Lemma 2.3 and Lemma 4.2, it follows
that $\mathcal A(k)  =\mathcal R$. \quad Q.E.D

%By Lemma 2.1 in [6] and the fact that $x \in \R_2^{\omega}$, we have
%that $0 =||xw-wx||_2 \ge ||(x-E_{(\mathcal A\otimes
%\B)^{\omega}}(x))w||_2 = ||(x-E_{(\mathcal A\otimes
%\B)^{\omega}}(x)) ||_2>0 $, which is a contradiction. Therefore
%$\mathcal R_2=0 $ and $\mathcal R=\mathcal R_1$. From Lemma 2.3 and
%4.3, it follows that $\mathcal A\otimes \B=\mathcal R$. Hence
%$\A\otimes \B$ is maximal injective in $\M(k)\otimes \B$.

\pb The following is the main result in this section.
\begin{Th}
Following the notations as above. Suppose  $\{n_i\}_{i=1}^N$ is a
sequence of integers where $n_i\ge 2$ for all $1\le i\le N$ and $N$
is finite or infinite. Let $F(n_i)$ be the free group with the
standard generators $\{g_{i,j}\}_{j=1}^{n_i}$ for all $1\le i\le N$.
Let the group $G $ be $\times_{i=1}^N  F(n_i)$, the direct product
of $F(n_1), \ldots, F(n_N)$. And $F(n_i)$ is identified with its
canonical image in $G$. Let $\lambda$ be the left regular
representation of $G$ and $\M=L(G)\cong \otimes_{i=1}^N L(F(n_i))$
be the group von Neumann algebra associated with $G$. Let $\A$ be
the abelian von Neumann subalgebra of $\M$ generated by the unitary
elements $\{\lambda(g_{i,1})| 1\le i\le N\}$.  Then $\A$ is a
maximal injective subalgebra of $\M$ and not contained in any
hyperfinite subfactor of $\M$.
\end{Th}

 The proof of Theorem 4.1 is divided into two
different cases: (i) $N$ is finite (ii) $N$ is infinite. Therefore
the theorem will follow easily from the following two propositions.
\begin{Prop} When $N$ is finite, $\A$ is a maximal injective
subalgebra of $\M$. Moreover, if $\B\cong L(\Bbb Z)$, then
$\A\otimes \B$ is a maximal injective abelian von Neumann subalgebra
of $\M\otimes \B$.
\end{Prop}

\begin{Prop} When $N$ is infinite, $\A$ is a maximal injective
subalgebra of $\M$.
\end{Prop}

 \noindent {\bf Proof of Proposition 4.1: }The proposition follows easily
from Lemma  4.4, Lemma 4.5, Lemma 4.6 and Lemma 4.7.

\vspace{0.5cm}

\noindent{\bf Proof of Proposition 4.2: } Recall $$\begin{aligned}
% H&= \text {subgroup of $G$ generated by }\{ g_{i,1}\}_{i=1}^\infty\\
H_i &=\text {subgroup of $F(n_i)$ generated by }  g_{i,1}, \qquad \text { for }  i\ge 1 \\
 %G_i& = \left (\times_{k=1}^{i-1}F(n_k)\right ) \times H_i \times \left
% (\times_{k=i+1}^\infty
% F_{n_k}\right),  \qquad \text { for } i\ge 1 \\
 J_i &=\left (\times_{k=1}^{i-1}F(n_k)\right )   \times \left
 (\times_{k=i}^\infty
 H_{k}\right),  \qquad \text { for }  i\ge 1 \\
\end{aligned}$$ %Therefore $\A=L(H)$.
  We denote
by  $\A_i$, or $\B_i$, the von Neumann subalgebra $L(H_i)$,  or
$L(J_i)$ respectively, of $\M$, for  $i\ge 1$ . \vspace{0.3cm}

\noindent It is easy to see that $ \B_i \nearrow \M$. Assume that
$\R$ is an injective von Neumann subalgebra of $\M$ that contains
$\A$ properly. Hence there exists some $x$ in $\R$ but not in $\A$.
There exists some positive number $a$ such that $||x||_2  > a
>||E_{\A}(x)||_2.$ Note that $E_{\B_i}(x) \rightarrow x$ when $i$
goes to infinity. There is some $k\in \Bbb N$ such that
$||E_{B_k}(x)||_2>a.$ Since $$\B_k\cong L(F(n_1))\otimes \cdots
\otimes L(F(n_{k-1}))\otimes \left ( \otimes_{j=k}^\infty A_j \right
),$$ by same arguments as in Lemma 4.7 we know that $E_{\B_k}(x)
=\frac 1 {2n}\lim_{l=-n}^n u^{l}xu^{-l}, $ where $u$ is the Haar
unitary element that generates $\otimes_{j=k}^\infty A_j .$ It
follows that $E_{\B_k}(x) \in \R\cap \B_k$. Denote $E_{\B_k}(x)$ by
$y$. Note that $||y||_2>a>||E_{\A}(y)||_2$. We know that $y$ is not
contained in $\A$. Let $\mathcal S$ be the von Neumann subalgebra of
$\R\cap \B_k$ generated by $y$ and $\A$. Since $\R$ is injective,
$\mathcal S$ is also injective. By Proposition 4.1, $\mathcal A$
($=\otimes_{i=1}^\infty \A_i$) is maximal injective in $\B_k \cong
L(F(n_1))\otimes \cdots \otimes L(F(n_{k-1}))\otimes \left (
\otimes_{j=k}^\infty A_j \right ).$ It contradicts with the fact
that $0\ne y \in \mathcal S\setminus \A$ and $\mathcal S$ is
injective. Hence $\A$ is a maximal injective von Neumann subalgebra
of $\M$.  Q.E.D

\pb {\bf Remark: } When $N$ is infinite, we obtain   examples of
McDuff factors of type II$_1$,   infinite tensor products of free
group factors,  that contains an abelian von Neumann subalgebra as
the maximal injective abelian von Neumann subalgebra. These McDuff
factors have self-adjoint operators that are not contained in any
hyperfinite subfactors, which also answers Kadison's problem \#7 in
the negative.

%\pb {\bf Remark: }  Let $\{n_i\}$ be a sequence of positive integers
%where $n_i\ge 2$ and $N$ is finite or infinite. Let $L(F(n_i))$ be
%the free group factor on $n_i$ generators. Let $\{\R_i\}$ is a
%family of non-atomic injective von Neumann algebras. Let
%$\M_i=L(F(n_i)) * \R_i$ for all $i$. Let $\M$ be the   of the tensor
%product of all von Neumann algebrash $\M_i$, .i.e.
%$\M=\otimes_{i=1}^N\M_i$, and $\R=\otimes_{i=1}^N \R_i$. Combining
%with the calculations in [3], Theorem 4.1 can be easily extended as
%follows.
%\begin{Col}
%With the notations as above, $\R$ is a maximal injective von Neumann
%subalgebra of $\M$.
%\end{Col}

%\vspace{3cm}
\vspace{2cm}
 \newpage
\noindent {\bf REFERENCES}

\begin{enumerate} [1.]

  \item A. Connes, ``Classification of Injective Factors,"
Ann. of Math. 104 (1976), 73-116.

%\item K. Dykema, ``Free Products of Hyperfinite von Neumann Algebras
%and Free Dimension," Duke Math. J. 69 (1993), 97-119.

\item L. Ge, ``On Maximal Injective Subalgebras of Factors,"
Advances in Mathematics 118 (1996), 34-70.

\item R. Kadison, ``Problems on  von  Neumann  Algebras,"  Baton
Rouge  Conference,  unpublished.

\item R. Kadison, ``Diagonalizing Matrices,"
 American Journal of Mathematics Vol. 106,
No. 6 (Dec., 1984), pp. 1451-1468."

\item R. Kadison, ``Fundamental of the Theory of Operator Algebras,"
Vols. 1 and 2, Academic Press, Orlando, 1983 and 1986.

\item D. McDuff, ``Central Sequences and the Hyperfinite Factors,''
Proc. London Math. Soc. 21 (1970), 443-461.

\item F. Murray and von Neumann, ``Rings of Operators,'' Ann. of
Math. 44 (1943), 716-808.

\item S.Popa, ``Orthogonal Pairs of subalgebras in Finite von
Neumann Algebras," J. Operator Theory 9 (1983), 253-268.

\item S. Popa, ``Maximal  Injective  Subalgebras  in Factors
Associated with Free groups," Advances in Mathematics 50, (1983),
27-48.

 \item S. Sakai,  ``$C^*$-algebras and $W^*$-algebras,''
  Reprint of the 1971 edition. Classics in Mathematics. Springer-Verlag, Berlin, 1998.

\item M. Takesaki, ``Theory of operator algebras," Vols. 1, 2 and 3,
 Encyclopaedia of Mathematical Sciences, 124, 125 and 127.  Springer-Verlag, Berlin, 2002, 2003 and 2003.

\end{enumerate}

\end{document}